\documentclass[11pt]{amsart}

\usepackage{amsfonts, amsmath, amssymb, amsthm}
\theoremstyle{plain}
\newtheorem{thm}{Theorem}[section]
\newtheorem{lem}[thm]{Lemma}
\newtheorem{prop}[thm]{Proposition}
\newtheorem{cor}[thm]{Corollary}

\theoremstyle{definition}

\newtheorem{exmp}[thm]{Example}

\theoremstyle{remark}
\newtheorem{rem}[thm]{Remark}

\usepackage{txfonts}
\usepackage{mathbbol}
\usepackage{ae}
\usepackage{color}
\usepackage{graphicx}
\usepackage{url}
\usepackage[abs]{overpic}
\usepackage{subfigure}
\usepackage[letterpaper]{anysize}
\usepackage[vlined]{algorithm2e}

\SetKwInOut{Input}{Input}
\SetKwInOut{Output}{Output}
\newcommand\frA{{\mathfrak A}}
\newcommand\cH{{\mathcal H}}
\newcommand\cK{{\mathcal K}}
\newcommand\cZ{{\mathcal Z}}
\newcommand\NN{{\mathbb N}}
\newcommand\RR{{\mathbb R}}
\newcommand\tropPG{{{\mathbb T}{\mathbb P}}}
\newcommand\tdet{\operatorname{tdet}}
\newcommand\tsgn{\operatorname{tsgn}}
\newcommand\tconv{\operatorname{tconv}}
\renewcommand\Vert{\operatorname{Vert}}
\newcommand\tr{{\operatorname{tr}}}
\newcommand\lr{\operatorname{bottom}}
\newcommand\rh{\operatorname{right}}
\newcommand\hl{\operatorname{top}}
\newcommand\lh{\operatorname{left}}
\newcommand\Sym{\operatorname{Sym}}
\newcommand\abs[1]{|#1|}
\newcommand\norm[1]{|\hskip-.2ex|#1|\hskip-.2ex|}
\newcommand\SetOf[2]{\left\{#1\,\vphantom{#2}\right.\left|\vphantom{#1}\,#2\right\}}

\begin{document}

\title{Tropical Halfspaces}
\author[Joswig]{Michael Joswig}
\address{Michael Joswig, Institut f\"ur Mathematik, Ma 6-2, TU Berlin,
  10623 Berlin, Germany}
\thanks{This work has been carried out while visiting the Mathematical Sciences Research Institute
    in Berkeley for the special semester on Discrete and Computational Geometry}
\email{joswig@math.tu-berlin.de}
\date{\today}

\begin{abstract}
  As a new concept tropical halfspaces are introduced to the (linear algebraic) geometry of the tropical semiring
  $(\RR,\min,+)$.  This yields exterior descriptions of the tropical polytopes that were recently studied by Develin and
  Sturmfels~\cite{TropConvex} in a variety of contexts.  The key tool to the understanding is a newly defined sign of
  the tropical determinant, which shares remarkably many properties with the ordinary sign of the determinant of a
  matrix.  The methods are used to obtain an optimal tropical convex hull algorithm in two dimensions.
\end{abstract}

\maketitle

\section{Introduction}

The set $\RR$ of real numbers carries the structure of a semiring if equipped with the \emph{tropical addition}
$\lambda\oplus\mu=\min\{\lambda,\mu\}$ and the \emph{tropical multiplication} $\lambda\odot\mu=\lambda+\mu$, where $+$
is the ordinary addition.  We call the triplet $(\RR,\oplus,\odot)$ the \emph{tropical semiring}\footnote{Other authors
  reserve the name \emph{tropical semiring} for $(\NN\cup\{+\infty\},\min,+)$ and call $(\RR\cup\{+\infty\},\min,+)$ the
  min-plus-semiring.}.  It is an equally simple and important fact that the operations $\oplus,\odot:\RR\times\RR\to\RR$
are continuous with respect to the standard topology of~$\RR$.  So the tropical semiring is, in fact, a topological
semiring.  Considering the tropical scalar multiplication
\[\lambda\odot (\mu_0,\dots,\mu_d)=(\lambda+\mu_0,\dots,\lambda+\mu_d)\] (and component-wise tropical addition) turns the
set $\RR^{d+1}$ into a semimodule.

The study of the linear algebra of the tropical semiring and, more generally, of idempotent semirings, has a long
tradition.  Applications to combinatorial optimization, discrete event systems, functional analysis etc.\ abound.  For an
introduction the reader is referred to the monograph by Baccelli et al.~\cite{MR94b:93001}.  A recent contribution in
the same vein, with many more references, is Cohen, Gaubert, and Quadrat~\cite{math.FA/0212294}.

Convexity in the tropical world (and even in a more general setting) was first studied by Zimmermann~\cite{MR56:11869}.
Following the approach of Develin and Sturmfels~\cite{TropConvex} here we stress the point of view of discrete geometry.
We recall some of the key definitions.  A subset $S\subset\RR^{d+1}$ is \emph{tropically convex} if for any two points
$x,y\in S$ the \emph{tropical line segment}
\[[x,y]=\SetOf{\lambda\odot x\oplus\mu\odot y}{\lambda,\mu\in\RR}\]
is contained in~$S$.  The \emph{tropical convex hull} of a set $S\subset\RR^{d+1}$ is the smallest tropically convex set
containing~$S$; it is denoted by $\tconv S$.  It is easy to see, cf.\ \cite[Proposition~4]{TropConvex}, that $\tconv
S=\SetOf{\lambda_1\odot x_1\oplus\dots\oplus\lambda_n\odot x_n}{\lambda_i\in\RR,\ x_i\in S}$.  A \emph{tropical
  polytope} is the tropical convex hull of finitely many points.  Since any convex set in $\RR^{d+1}$ is closed under
tropical multiplication with an arbitrary scalar, it is common to identify tropically convex sets with their respective
images under the canonical projection onto the $d$-dimensional \emph{tropical projective space}
\[\tropPG^d=\SetOf{\RR\odot x}{x\in\RR^{d+1}}=\RR^{d+1}/\RR(1,\dots,1).\] In explicit computations we often choose
\emph{canonical coordinates} for a point $x\in\tropPG^d$, meaning the unique non-negative vector in the class $\RR\odot
x$ which has at least one zero coordinate.  For visualization purposes, however, we usually normalize the coordinates by
choosing the first one to be zero (which can then be omitted): This identification
$(\xi_0,\dots,\xi_d)\mapsto(\xi_1-\xi_0,\dots,\xi_d-\xi_0):\tropPG^d\to\RR^d$ is a homeomorphism.

Develin and Sturmfels observed that \emph{tropical simplices}, that is tropical convex hulls of $d+1$ points
in~$\tropPG^d$ (in sufficiently general position), are related to Isbell's~\cite{MR32:431} \emph{injective envelope} of
a finite metric space, cf.\ \cite[Theorem~29]{TropConvex} and the Erratum.  Isbell's injective envelope in turn
coincides with the \emph{tight span} of a finite metric space that arose in the work of Dress and others; see the
paper~\cite{MR2003g:54077} and its list of references.  In a way, tropical simplices may be understood as non-symmetric
analogues of injective hulls or tight spans.

Additionally, tropical polytopes are interesting also from a purely combinatorial point of view: They bijectively
correspond to the regular polyhedral subdivisions of products of simplices; see \cite[Theorem 1]{TropConvex}.

The present paper studies tropical polytopes as geometric objects in their own right.  It is shown that, at least to
some extent, it is possible to develop a theory of tropical polytopes in a fashion similar to the theory of ordinary
convex polytopes.  The key concept introduced to this end is the notion of a tropical halfspace.  One of our main
results, Theorem~\ref{thm:halfspace}, gives a characterization of tropical halfspaces in terms of the tropical
determinant, which is the same as the min-plus-permanent already studied by Yoeli~\cite{MR23:A3768} and others; see also
Sturmfels, Richter-Gebert, and Theobald~\cite{FirstSteps}.  The proof leads to the definition of the \emph{faces} of a
tropical polytope in a natural way.  In the investigation, in particular, we prove that the faces form a distributive
lattice, cf.\ Theorem~\ref{thm:lattice}.  Moreover, as one would expect by analogy to ordinary convex polytopes, the
tropical polytopes are precisely the bounded intersections of finitely many tropical halfspaces, cf.\ 
Theorem~\ref{thm:main}.

It is a further consequence of our results on tropical polytopes that some concepts and ideas from computational
geometry can be carried over from ordinary convex polytopes to tropical polytopes.  In Section~\ref{sec:2d} this leads
us to a comprehensive solution of the convex hull problem in~$\tropPG^2$.  The general tropical convex hull problem in
arbitrary dimension is certainly interesting, but this is beyond our current scope.

The paper closes with a selection of open questions.

\section{Hyperplanes and Halfspaces}

We start this section with some observations concerning the topological aspects of tropical convexity.  As already
mentioned the tropical projective space $\tropPG^d$ is homeomorphic to $\RR^d$ with the usual topology.  Moreover, the
space $\tropPG^d$ carries a natural metric: For a point $x\in\tropPG^d$ with canonical coordinates $(\xi_0,\dots,\xi_d)$
let \[\norm{x}=\max\{\xi_0,\dots,\xi_d\}\] be the \emph{tropical norm} of~$x$.  Equivalently, for arbitrary coordinates
$(\xi_0',\dots,\xi_d')\in\RR\odot x$ we have that $\norm{x}=\max\SetOf{\abs{\xi_i'-\xi_j'}}{i\ne j}$.  We prove a
special case of~\cite[Theorem~17]{math.FA/0212294}:

\begin{lem}
  The map \[\tropPG^d\times\tropPG^d\to\RR:(x,z)\mapsto\norm{x-z}\] is a metric.
\end{lem}

\begin{proof}
  By definition the map is non-negative.  Moreover, it is clearly definite and symmetric.  We prove the triangle
  inequality: Assume that $x=(\xi_0,\dots,\xi_d)$, $z=(\zeta_0,\dots,\zeta_d)$, and that $y=(\eta_0,\dots,\eta_d)$ be a
  third point.  Then
  \begin{eqnarray*}
    \norm{x-z}&=&\max\SetOf{\abs{(\xi_i-\zeta_i)-(\xi_j-\zeta_j)}}{i\ne j}\\
    &=&\max\SetOf{\abs{(\xi_i-\xi_j)-(\eta_i-\eta_j)+(\eta_i-\eta_j)-(\zeta_i-\zeta_j)}}{i\ne j}\\
    &\le& \max\SetOf{\abs{(\xi_i-\xi_j)-(\eta_i-\eta_j)}+\abs{(\eta_i-\eta_j)-(\zeta_i-\zeta_j)}}{i\ne j}\\
    &\le&\max\SetOf{\abs{(\xi_i-\eta_i)-(\xi_j-\eta_j)}}{i\ne j} +
    \max\SetOf{\abs{(\eta_i-\zeta_i)-(\eta_j-\zeta_j)}}{i\ne j}\\
    &=&\norm{x-y}+\norm{y-z}.
  \end{eqnarray*}
\end{proof}

The topology induced by this metric coincides with the quotient topology on~$\tropPG^d$ (and thus with the natural
topology of~$\RR^d$).  In particular, $\tropPG^d$ is locally compact and a set $C\subset\tropPG^d$ is compact if and
only if it is closed and bounded.  Tacitly we will always assume that $d\ge2$.

\begin{prop}\label{prop:closure}
  The topological closure of a tropically convex set is tropically convex.
\end{prop}

\begin{proof}
  Let $S$ be a tropically convex set with closure $\bar S$.  Then, by \cite[Proposition~4]{TropConvex}, $\tconv(\bar S)$
  is the set of points in~$\tropPG^d$ which can be obtained as tropical linear combinations of points in~$\bar S$.  Now
  the claim follows from the fact that tropical addition and multiplication are continuous.
\end{proof}

From the named paper by Develin and Sturmfels we quote a few results which will be useful in our investigation.

\begin{thm}{\rm (\cite[Theorem 15]{TropConvex})}\label{thm:canonical-decomposition}
  A tropical polytope has a canonical decomposition as a finite ordinary polytopal complex, where the cells are both
  ordinary and tropical polytopes.
\end{thm}

\begin{prop}{\rm (\cite[Proposition 20]{TropConvex})}\label{prop:polytope-intersection}
  The intersection of two tropical polytopes is again a tropical polytope.
\end{prop}

\begin{prop}{\rm (\cite[Proposition 21]{TropConvex}; see also Helbig~\cite{MR89g:46118})}\label{prop:vertices}
  For each tropical polytope $P$ there is a unique minimal set $\Vert(P)\subset P$ with $\tconv(\Vert(P))=P$.
\end{prop}

The elements of $\Vert(P)$ are called the \emph{vertices} of~$P$.  The following is implied by
Theorem~\ref{thm:canonical-decomposition}.  There is also an easy direct proof which we omit, however.

\begin{prop}
  A tropical polytope is compact.
\end{prop}

The \emph{tropical hyperplane} defined by the \emph{tropical linear form} $a=(\alpha_0,\dots,\alpha_d)\in\RR^{d+1}$ is
the set of points $(\xi_0,\dots,\xi_d)\in\tropPG^d$ such that the minimum
\[
  \min\{\alpha_0+\xi_0,\dots,\alpha_d+\xi_d\}=\alpha_0\odot\xi_0\oplus\dots\oplus\alpha_d\odot\xi_d
\]
is attained at least twice.  The point $-a$ is contained in the tropical hyperplane defined by~$a$, and it is called its
\emph{apex}.  Note that any two tropical hyperplanes only differ by a translation.

\begin{prop}{\rm (\cite[Proposition 6]{TropConvex})}
  Tropical hyperplanes are tropically convex.
\end{prop}

The complement of a tropical hyperplane $\cH$ in $\tropPG^d$ has $d+1$ connected components corresponding to the facets
of an ordinary $d$-simplex.  We call each such connected component an \emph{open sector} of~$\cH$.  The topological
closure of an open sector is a \emph{closed sector}.  It is easy to prove that each (open or closed) sector is convex
both in the ordinary and in the tropical sense.

\begin{exmp}\label{exmp:sectors}
  Consider the zero tropical linear form $0\in\RR^{d+1}$.  The open sectors of the corresponding tropical hyperplane
  $\cZ$ are the sets $S_0,\dots,S_d$, where
  \[S_i=\SetOf{(\xi_0,\dots,\xi_d)}{\text{$\xi_i<\xi_j$ for all $j\ne i$}}.\]
  The closed sectors are the sets $\bar
  S_0,\dots,\bar S_d$, where \[\bar S_i=\SetOf{(\xi_0,\dots,\xi_d)}{\text{$\xi_i\le\xi_j$ for all $j\ne i$}}.\]
  In
  canonical coordinates this can be expressed as follows: \[S_i=\SetOf{(\xi_0,\dots,\xi_d)}{\text{$\xi_i=0$ and
      $\xi_j>0$, for $j\ne i$}}\]
  and \[\bar S_i=\SetOf{(\xi_0,\dots,\xi_d)}{\text{$\xi_i=0$ and $\xi_j\ge0$, for $j\ne
      i$}}.\]
\end{exmp}

Like any two tropical hyperplanes are related by a translation, each translation of a sector is again a sector.  We call
such sectors \emph{parallel}.

The following simple observation is one of the keys to the structural results on tropical polytopes in the subsequent
sections.  It characterizes the solvability of one tropical linear equation.  For related results see Akian, Gaubert,
and Kolokoltsov~\cite{math.FA/0403441}.

\begin{prop}\label{prop:zero}
  Let $x_1,\dots,x_n\in\tropPG^d$.  Then $0\in\tconv\{x_1,\dots,x_n\}$ if and only if each closed sector $\bar S_k$ of
  the zero tropical linear form contains at least one $x_i$.
\end{prop}

\begin{proof}
  We write $\xi_{ij}$ for the canonical coordinates of the $x_i$ in~$\RR^{d+1}$.  Then all the $n(d+1)$ entries in the matrix
  \[
  \begin{pmatrix}
    \xi_{10} & \cdots & \xi_{1d}\\ \vdots & \ddots & \vdots\\ \xi_{n0} & \cdots & \xi_{nd}
  \end{pmatrix}
  \]
  are non-negative.  Hence \[0=\lambda_1\odot x_1\oplus\dots\oplus\lambda_n\odot x_n\] (with $\lambda_i\ge 0$, as we may
  assume without loss of generality) if and only if $\min\{\lambda_1+\xi_{1k},\dots,\lambda_n+\xi_{nk}\}=0$ for all $k$.
  We conclude that zero is in the tropical convex hull of $x_1,\dots,x_n$ if and only if for all $k$ there is an $i$
  such that $\xi_{ik}=0$ or, equivalently, $x_i\in\bar S_k$.
\end{proof}

Throughout the following we abbreviate $[d+1]=\{0,\dots,d\}$, and we write $\Sym(d+1)$ for the symmetric group of degree
$d+1$ acting on the set $[d+1]$.  Let $e_i$ be the $i$-th unit vector of~$\RR^{d+1}$.  Observe that under the natural
action of $\Sym(d+1)$ on~$\tropPG^d$ by permuting the unit vectors tropically convex sets get mapped to tropically
convex sets.  The set of all $k$-element subsets of a set~$\Omega$ is denoted by $\binom{\Omega}{k}$.

We continue our investigation with the construction of a two-parameter family of tropical polytopes.

\begin{exmp}
  We define the $k$-th \emph{tropical hypersimplex} in $\tropPG^d$ as
  \[\Delta_k^d=\tconv\SetOf{\sum_{i\in J}-e_i}{J\in\binom{[d+1]}{k}}\subset\tropPG^d.\]
  It is essential that
  \[\Vert(\Delta_k^d)=\SetOf{\sum_{i\in J}-e_i}{J\in\binom{[d+1]}{k}},\]
  for all $k>0$: This has to do with the fact
  that the symmetric group $\Sym(d+1)$ acts on the set, due to which either all or none of the points $\sum_{i\in
    J}-e_i$ is a vertex.  But from Proposition~\ref{prop:vertices} we know that
  $\emptyset\ne\Vert(\Delta_k^d)\subseteq\SetOf{\sum_{i\in J}-e_i}{J\in\binom{[d+1]}{k}}$, and hence the claim follows.
  Develin, Santos, and Sturmfels~\cite{TropRank} construct tropical polytopes from matroids.  The tropical
  hypersimplices arise as the special case of uniform matroids.
\end{exmp}

It is worth-while to look at two special cases of the previous construction.

\begin{exmp}\label{exmp:standard}
  The first tropical hypersimplex in $\tropPG^d$ is the $d$-dimensional \emph{tropical standard simplex}
  $\Delta^d=\Delta_1^d=\tconv\{-e_0,\dots,-e_d\}$.  Note that $\Delta^d$ is a tropical polytope which at the same time
  is an ordinary polytope.
\end{exmp}

\begin{exmp}\label{exmp:Ed}
  The second tropical hypersimplex $\Delta_2^d\subset\tropPG^d$ is the tropical convex hull of the $\binom{d+1}{2}$
  vectors $-e_i-e_j$ for all pairs $i\ne j$.  The tropical polytope $\Delta_2^d$ is not convex in the ordinary sense.
  It is contained in the tropical hyperplane $\cZ$ corresponding to the zero tropical linear form.  For $d=2$ see
  Figure~\ref{fig:full-triangle} below.
\end{exmp}

\begin{prop}\label{prop:Ed}
  The second tropical hypersimplex $\Delta_2^d\subset\tropPG^d$ is the intersection of the tropical hyperplane $\cZ$
  corresponding to the zero tropical linear form with the set of points whose tropical norm is bounded by~$1$.
\end{prop}

\begin{proof}
  Clearly, $-e_i-e_j\in\cZ$ for $i\ne j$.  We have to show that a point $x$ with canonical coordinates
  $(\xi_0,\dots,\xi_d)$ and $\norm{x}\le 1$ such that, e.g., $\xi_0=\xi_1=0$, is a tropical linear combination of the
  $\binom{d+1}{2}$ vertices of~$\Delta_2^d$.  We compute
  \[x=(0,0,1,\dots,1)\oplus\xi_2\odot(0,1,0,1,\dots,1)\oplus\dots\oplus\xi_d\odot(0,1,\dots,1,0),\] and hence the claim.
\end{proof}

In particular, this implies that $\Delta_2^d$ contains $\Delta_k^d$, for all $k>2$.  A similar computation further shows
that $\Delta_k^d\supsetneq\Delta_{k+1}^d$, for all~$k$.

\begin{exmp}
  The ordinary $d$-dimensional $\pm1$-cube \[C^d=\SetOf{(0,\xi_1,\dots,\xi_d)}{-1\le\xi_i\le 1}\]
  is a tropical
  polytope: $C^d=\tconv\{-e_0-2e_1,\dots,-e_0-2e_d,e_1+\dots+e_d\}$.
\end{exmp}

One way of reading Proposition~\ref{prop:Ed} is that the intersection of the tropical hyperplane corresponding to the
zero tropical linear form with the ordinary $\pm1$-cube is a tropical polytope.  An important consequence is the
following.

\begin{cor}\label{cor:hyperplane-intersection}
  The (non-empty) intersection of a tropical polytope with a tropical hyperplane is again a tropical polytope.
\end{cor}

\begin{proof}
  Let $P\subset\tropPG^d$ be a tropical polytope and $\cH$ a tropical hyperplane.  As usual, up to a translation we can
  assume that $\cH=\cZ$ corresponds to the zero tropical linear form.  By Proposition~\ref{prop:Ed} the intersection
  $P\cap\cZ$ is contained in a suitably scaled copy of the second tropical hypersimplex~$\Delta_2^d$.  Now the claim
  follows from Proposition~\ref{prop:polytope-intersection}.
\end{proof}

A \emph{closed tropical halfspace} in~$\tropPG^d$ is the union of at least one and at most $d$ closed sectors of a fixed
tropical hyperplane.  Hence it makes sense to talk about the \emph{apex} of a tropical halfspace.  An \emph{open
  tropical halfspace} is the complement of a closed one.  Clearly, the topological closure of an open tropical halfspace
is a closed tropical halfspace.  To each (open or closed) tropical halfspace $\cH^+$ there is an \emph{opposite} (open
or closed) tropical halfspace $\cH^-$ formed by the sectors of the corresponding tropical hyperplane which are not
contained in~$\cH^+$.  Two halfspaces are \emph{parallel} if they are formed of parallel sectors.

\begin{lem}\label{lem:move}
  Let $a+\bar S_k\subset\tropPG^d$ be a closed sector, for some $k\in[d+1]$, and $b\in a+\bar S_k$ a point inside.  Then
  the parallel sector $b+\bar S_k$ is contained in $a+\bar S_k$.
\end{lem}

Note that this includes the case where $b$ is a point in the boundary $a+(\bar S_k\setminus S_k)$.  The proof of the
lemma is omitted.

\begin{prop}
  Each closed tropical halfspace is tropically convex.
\end{prop}

\begin{proof}
  Let $\cH^+$ be a closed tropical halfspace.  Without loss of generality, we can assume that $\cH^+$ is the union the of
  closed sectors $\bar S_{i_1}, \dots, \bar S_{i_l}$ of the tropical hyperplane $\cZ$ corresponding to the zero tropical
  linear form.  Since we already know that each $\bar S_{i_k}$ is tropically convex, it suffices to consider, e.g.,
  $x\in\bar S_{i_1}$ and $y\in\bar S_{i_2}$ and to prove that $[x,y]\subset\cH^+$.  Let $(\xi_0,\ldots,\xi_d)$ and
  $(\eta_0,\dots,\eta_d)$ be the canonical coordinates of $x,y\in\tropPG^d$, respectively.  Since $x\in\bar S_{i_1}$ and
  $y\in\bar S_{i_2}$ we have that $\xi_{i_1}=0$ and $\eta_{i_2}=0$.  Then the minimum
  \[\min\{\lambda+\xi_0,\dots,\lambda+\xi_d,\mu+\eta_0,\dots,\mu+\eta_d\}\] is $\lambda=\lambda+\xi_{i_1}$ or
  $\mu=\mu+\eta_{i_2}$, for arbitrary $\lambda,\mu\in\RR$.  This is equivalent to $\lambda\odot x+\mu\odot y\in\bar
  S_{i_1}\cup\bar S_{i_2}$, which implies the claim.
\end{proof}

A similar argument shows that open tropical halfspaces are tropically convex.

\begin{cor}
  The boundary of a tropical halfspace is tropically convex.
\end{cor}

\begin{proof}
  The boundary of a closed tropical halfspace $\cH^+$ is the intersection of $\cH^+$ with its opposite closed tropical
  halfspace~$\cH^-$.
\end{proof}

\theoremstyle{plain}
\newtheorem{TST}[thm]{Tropical Separation Theorem}
\begin{TST}\label{thm:TST}
  Let $P$ be tropical polytope, and $x\not\in P$ a point outside.  Then there is a closed tropical halfspace
  containing $P$ but not~$x$.
\end{TST}

\begin{proof}
  From Proposition~\ref{prop:zero} we infer that there is a closed sector $x+\bar S_k$ of the tropical hyperplane with
  apex~$x$ which is disjoint from~$P$.  Now $e_k$ is the unique coordinate vector such that $e_k\not\in\bar S_k$.  Since
  $P$ is compact and $\bar S_k$ is closed there is some $\epsilon>0$ such that the closed sector $x+\epsilon e_k+\bar
  S_k$ is disjoint from~$P$.  The complement of the open sector $x+\epsilon e_k+S_k$ is a closed tropical halfspace of
  the desired kind.
\end{proof}

Tropical halfspaces implicitly occur in the work of Cohen, Gaubert, and Quadrat~\cite{math.FA/0212294}.  In particular,
their results imply the Tropical Separation Theorem.  In fact, a variation of this result already occurs in
Zimmermann~\cite{MR56:11869}.  Another variant of the same is the Tropical Farkas Lemma of Develin and
Sturmfels~\cite[Proposition~9]{TropConvex}.

\section{Exterior Descriptions of Tropical Polytopes}

Throughout this section let $P\subset\tropPG^d$ be a tropical polytope.  Like their ordinary counterparts tropical
polytopes also have an exterior description.

\begin{lem}\label{lem:first-outer-description}
  The tropical polytope $P$ is the intersection of the closed tropical halfspaces that it is contained in.
\end{lem}

\begin{proof}
  Let $P'$ be the intersection of all the tropical halfspaces which contain~$P$.  Clearly, $P'$ contains~$P$.  Suppose
  that there is a point $x\in P'\setminus P$.  Then, by the Tropical Separation Theorem, there is a closed tropical
  halfspace which contains~$P$ but not~$x$.  This contradicts the assumption that $P'$ is the intersection of all such
  tropical halfspaces.
\end{proof}

Of course, the set of closed tropical halfspaces that contain the given tropical polytope~$P$ is partially ordered by
inclusion.  A closed tropical halfspace is said to be \emph{minimal} with respect to~$P$ if it is a minimal element in
this partial order.

A key observation in what follows is that the minimal tropical halfspaces come from a small set of candidates only.  For
a given finite set of points $p_1,\dots,p_n\in\tropPG^d$ let the \emph{standard affine hyperplane arrangement} be
generated by the ordinary affine hyperplanes \[p_i+\SetOf{(0,\xi_1,\dots,\xi_d)\in\RR^{d+1}}{\xi_j=0}\ \text{and}\ 
p_i+\SetOf{(0,\xi_1,\dots,\xi_d)\in\RR^{d+1}}{\xi_j=\xi_k}.\] For an example illustration see
Figure~\ref{fig:arrangement}.  A \emph{pseudovertex} of~$P$ is a vertex of the standard affine hyperplane arrangement
with respect to $\Vert(P)$ which is contained in the boundary~$\partial P$.  In~\cite{TropConvex} our pseudovertices are
called the \emph{vertices}.

The following is a special case of~\cite[Proposition 18]{TropConvex}.

\begin{prop}\label{prop:bounded-cells}
  The bounded cells of the standard affine hyperplane arrangement are tropical polytopes which are at the same time
  ordinary convex polytopes.
\end{prop}

\begin{prop}\label{prop:pseudo-vertex}
  The apex of a closed tropical halfspace that is minimal with respect to~$P$ is a pseudovertex of~$P$.
\end{prop}

\begin{proof}
  Let $\cH^+$ be a closed tropical halfspace, with apex $a$, which minimally contains~$P$.  Suppose that $a$ is not a
  vertex of the standard affine hyperplane arrangement~$\frA$ generated by~$\Vert(P)$, but rather $a$ is contained in
  the relative interior of some cell~$C$ of~$\frA$ of dimension at least one.  Now there is some $\epsilon>0$ such that
  for each point $a'$ in~$C$ with $\norm{a'-a}<\epsilon$ the closed tropical halfplane with apex~$a'$ and parallel
  to~$\cH^+$ still contains~$P$.  For each $a'\in\cH^+$ the corresponding translate is contained in $\cH^+$ and hence
  $\cH^+$ is not minimal.  Contradiction.
  
  It remains to show that $a\in P$.  Again suppose the contrary.  Then, by the Tropical Separation
  Theorem~\ref{thm:TST}, there is a closed halfspace $\cH_1^+$ containing~$P$ but not~$a$.  Now, since $\cH^+$ is
  minimal, $\cH_1^+$ is not contained in $\cH^+$ and, in particular, $\cH_1^+$ is not parallel to $\cH^+$.  As
  $a\not\in\cH_1^+$ the closed tropical halfspace $\cH_2^+$ with apex~$a$ which is parallel to $\cH_1^+$ contains~$P$.
  We infer that $\cH^+\cap\cH_2^+\subsetneq\cH^+$ is a closed tropical halfspace (with apex $a$) which contains~$P$.
  This contradicts the minimality of~$\cH^+$.
\end{proof}

\begin{cor}\label{cor:finite}
  There are only finitely many closed tropical halfspace which are minimal with respect to~$P$.
\end{cor}

\begin{proof}
  The standard affine hyperplane arrangement generated by $\Vert(P)$ is finite, and thus there are only finitely many
  pseudovertices.  Since there are only $2^{d+1}-2$ closed affine halfspaces with a given apex,\footnote{The
    Example~\ref{exmp:full-triangle} shows that there may indeed be more than one minimal halfspace with a given apex.}
  the claim now follows from Proposition~\ref{prop:pseudo-vertex}.
\end{proof}

This immediately gives the following.

\begin{cor}\label{cor:irredundant-outer-description}
  The tropical polytope $P$ is the intersection of the (finitely many) minimal closed tropical halfspaces that it is
  contained in.
\end{cor}

We can now prove our first main result.

\begin{thm}\label{thm:main}
  The tropical polytopes are exactly the bounded intersections of finitely many tropical halfspaces.
\end{thm}

\begin{proof}
  Let $P$ be the bounded intersection of finitely many tropical halfspaces $\cH_1^+,\dots,\cH_m^+$.  Then $P$ is the
  union of (finitely many) bounded cells of the standard affine hyperplane arrangement corresponding to the apices of
  $\cH_1^+,\dots,\cH_m^+$.  By Proposition~\ref{prop:bounded-cells} each of those cells is the tropical convex hull of
  its pseudovertices.  Since $P$ is tropically convex, this property extends to~$P$, and $P$ is a tropical polytope.

  The converse follows from Corollary~\ref{cor:irredundant-outer-description}.
\end{proof}

Ordinary polytope theory is combinatorial to a large extent.  This is due to the fact that many important properties of
an ordinary polytope are encoded into its face lattice.  While it is tempting to start a combinatorial theory of
tropical polytopes from the results that we obtained so far, this turns out to be quite intricate.  Here we give a brief
sketch, while a more detailed discussion will be picked up in a forthcoming second paper.

A \emph{boundary slice} of the tropical polytope~$P$ is the tropical convex hull of $\Vert(P)\cap\partial\cH^+$ where
$\cH^+$ is a closed tropical halfspace containing~$P$.  The boundary slices are partially ordered by inclusion; we call
a maximal element of this finite partially ordered set a \emph{facet} of~$P$.  Let $F_1,\dots,F_m$ be facets of~$P$.
Then the set
\[F_1\sqcap\dots\sqcap F_m=\tconv(\Vert(P)\cap F_1\dots\cap F_m)\] is called a \emph{proper face} of~$P$ provided that
$F_1\sqcap\dots\sqcap F_m\ne\emptyset$.  The sets $\emptyset$ and~$P$ are the \emph{non-proper faces}.  The faces of a
tropical polytope are partially ordered by inclusion, the maximal proper faces being the facets.  Note that, by
definition, faces of tropical polytopes are again tropical polytopes.

\begin{thm}\label{thm:lattice}
  The face poset of a tropical polytope is a finite distributive lattice.
\end{thm}

\begin{proof}
  We can extend the definition of~$\sqcap$ to arbitrary faces, this gives the \emph{meet} operation.  There is no choice
  for the \emph{join} operation then: $G\sqcup H$ is the meet of all facets containing $G$ and $H$, for arbitrary faces
  $G$ and~$H$.  Denote the set of facets containing the face~$G$ by $\Phi(G)$, that is, $G=\bigsqcap\Phi(G)$.  It is
  immediate from the definitions that $\Phi(G\sqcup H)=\Phi(G)\cap\Phi(H)$ and $\Phi(G\sqcap H)=\Phi(G)\cup\Phi(H)$.
  Hence the absorption and distributive laws are inherited from the boolean lattice of subsets of the set of all facets.
\end{proof}

In order to simplify some of the discussion below we shall introduce a certain non-degeneracy condition: A set
$S\subset\tropPG^d$ is called \emph{full}, if it is not contained in the boundary of any tropical halfspace.  If $S$ is
not contained in any tropical hyperplane, then, clearly, $S$ is full.  As the Example~\ref{exmp:full-triangle} below
shows, however, the converse does not hold.

\begin{figure}[htbp]
  \begin{center}
    \subfigure[Tropical standard simplex $\Delta^2$.\label{fig:delta2}]{%
      \includegraphics[width=.35\textwidth]{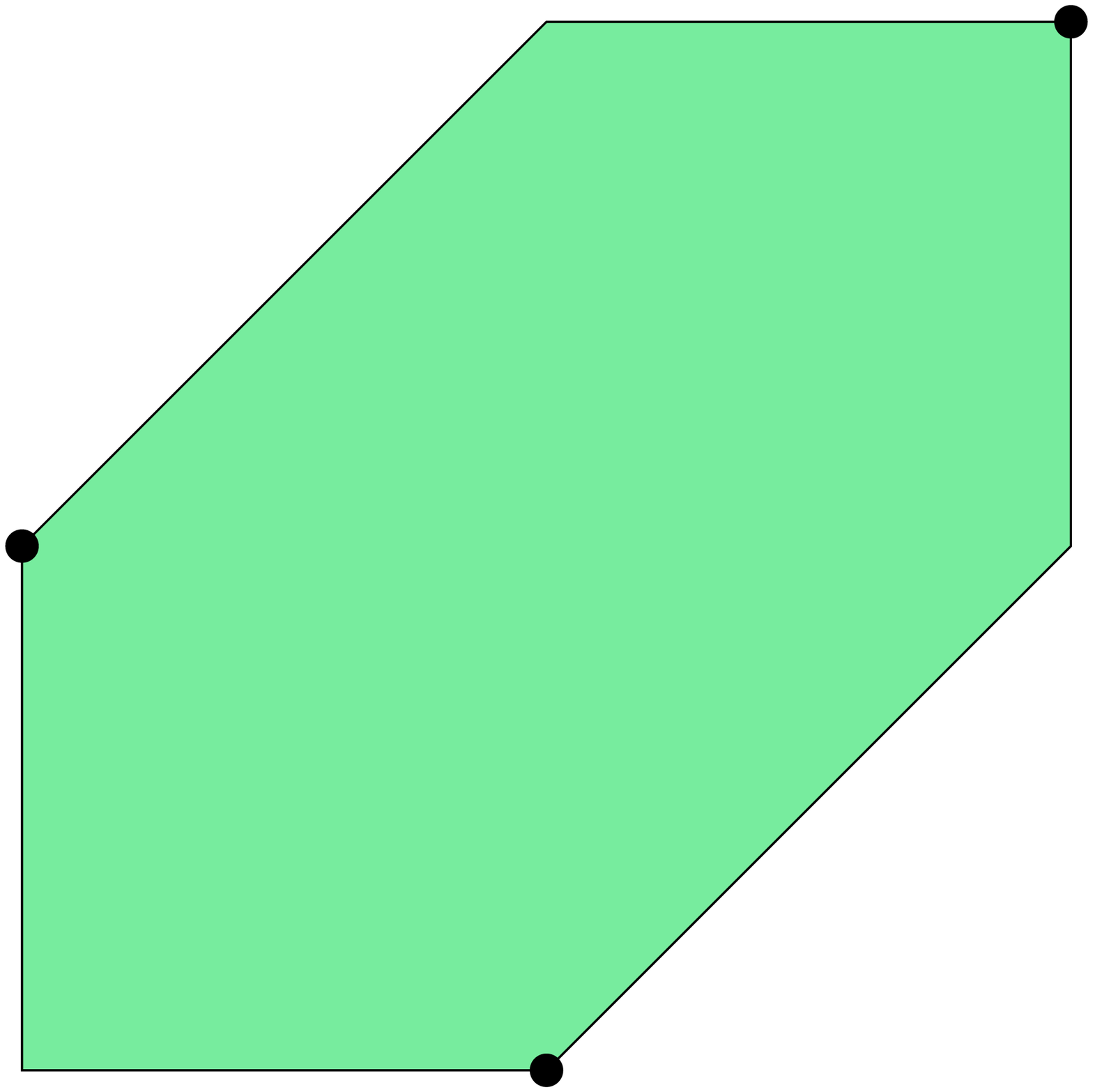}}
    \quad
    \subfigure[The second tropical hypersimplex~$\Delta_2^2\subset\tropPG^2$ has collinear vertices.\label{fig:full-triangle}]{%
      \includegraphics[width=.35\textwidth]{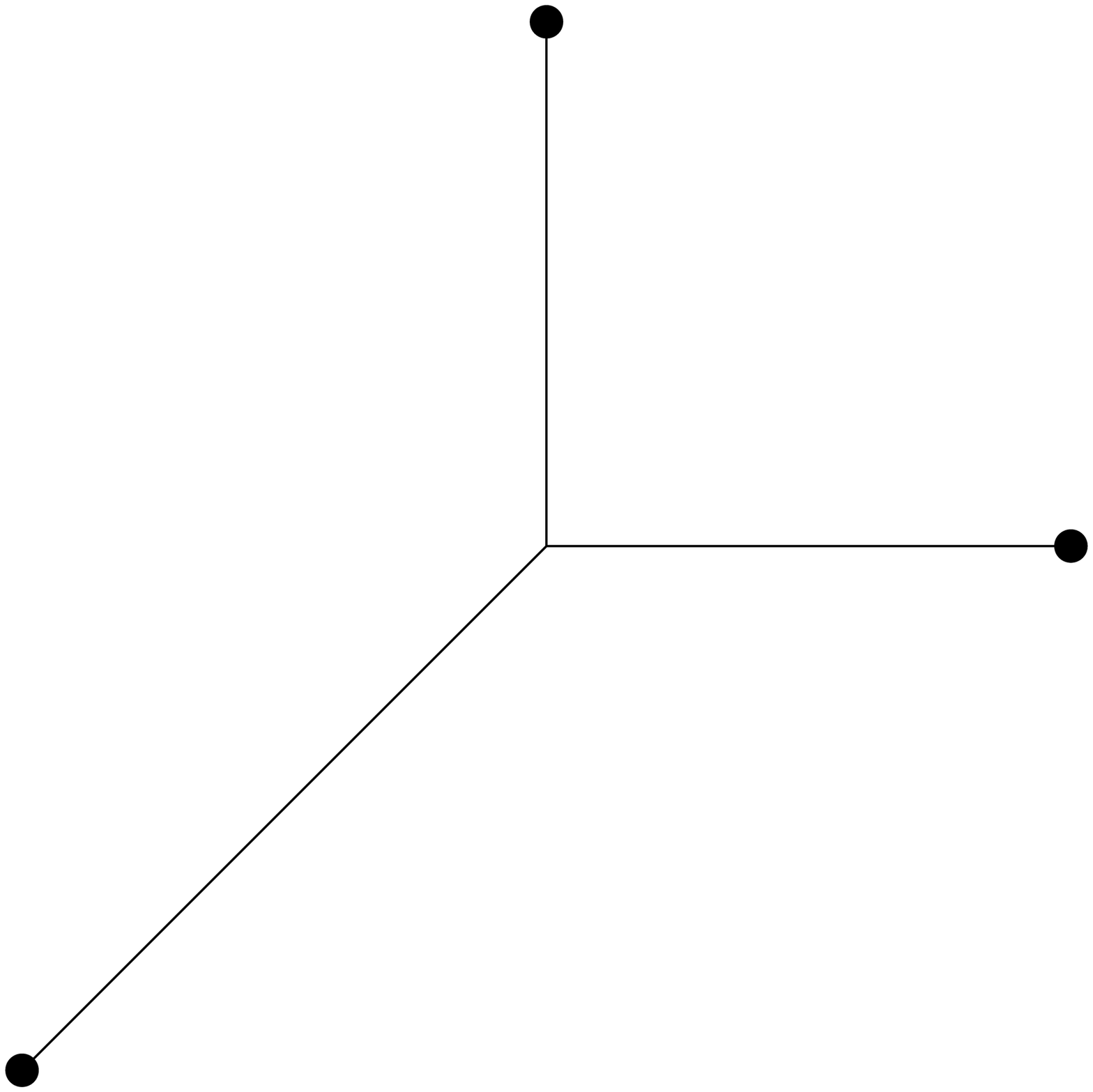}}
    \caption{Two full tropical triangles in $\tropPG^2$.  Both have the same face lattice as an ordinary triangle.}
  \end{center}
\end{figure}

\begin{exmp}
  The tropical standard simplex \[\Delta^2=\tconv\{(0,1,1),(1,0,1),(1,1,0)\}\] is not contained in a tropical
  hyperplane, and hence it is full; see Figure~\ref{fig:delta2}.  It is the intersection of the three minimal closed
  tropical halfspaces $(1,0,0)+\bar S_0$, $(0,1,0)+\bar S_1$, and $(0,0,1)+\bar S_2$.  The tropical line segments
  $[(0,1,1),(1,0,1)]$, $[(1,0,1),(1,1,0)]$, and $[(1,1,0),(0,1,1)]$ form the facets.  The three vertices form the only
  other proper faces.
\end{exmp}

\begin{exmp}\label{exmp:full-triangle}
  The second tropical hypersimplex \[\Delta_2^2=\tconv\{(1,0,0),(0,1,0),(0,0,1)\}\] is a full tropical triangle in the
  tropical plane~$\tropPG^2$, which is contained in the tropical line~$\cZ$ corresponding to the zero tropical linear
  form.  There are \emph{six} minimal closed tropical halfspaces: $\bar S_0\cup\bar S_1$, $\bar S_1\cup\bar S_2$, $\bar
  S_0\cup\bar S_2$, $(1,0,0)+\bar S_0$, $(0,1,0)+\bar S_1$, and $(0,0,1)+\bar S_2$.  Note that the three closed tropical
  halfspaces $\bar S_0\cup\bar S_1$, $\bar S_1\cup\bar S_2$, and $\bar S_0\cup\bar S_2$ share the origin as their apex.
  The tropical line segments $[(1,0,0),(0,1,0)]$, $[(0,1,0),(0,0,1)]$, and $[(0,0,1),(1,0,0)]$ form the facets.  Like in
  the example above the three vertices form the only other proper faces.  The triangle is depicted in
  Figure~\ref{fig:full-triangle}.
\end{exmp}

\begin{rem}  
  It is a consequence of Proposition~\ref{prop:polytope-intersection} and Corollary~\ref{cor:hyperplane-intersection}
  that boundary slices of tropical polytopes are tropical polytopes, but they need not be faces: E.g., the intersection
  of $\Delta_2^2$ with the boundary of the halfspace $(0,0,1/2)+(\bar S_1\cup\bar S_2)$ is the tropical (and at the same
  time ordinary) line segment $[(0,0,1/2),(0,0,1)]$ which contains the face $(0,0,1)$ and is properly contained in the
  intersection $[(1,0,0),(0,0,1)]\cap[(0,1,0),(0,0,1)]$ of two facets.
\end{rem}

\begin{rem}
  It is easy to see that the face lattice of a tropical $n$-gon in $\tropPG^2$ (which is necessarily full) is always the
  same as that of an ordinary $n$-gon.  This will play a role in the investigation of the algorithmic point of view in
  Section~\ref{sec:2d}.
\end{rem}

Minimal closed tropical halfspaces can be recognized by the intersection of their boundaries with~$P$, provided that $P$
is full.

\begin{prop}\label{prop:unique}
  Assuming that $P$ is full, let $\cH_1^+$ and $\cH_2^+$ be closed tropical halfspaces which are minimal with respect
  to~$P$.  If $\partial\cH_1^+\cap P=\partial\cH_2^+\cap P$ then $\cH_1^+=\cH_2^+$.
\end{prop}

\begin{proof}
  If $P$ is full then it is impossible that for any closed tropical halfspace $\cH^+$ containing $P$ the opposite closed
  halfspace~$\cH^-$ also contains~$P$.
  
  Let $a_1$ and $a_2$ be the respective apices of $\cH_1^+$ and $\cH_2^+$.  By Lemma~\ref{prop:pseudo-vertex} the points
  $a_1$ and $a_2$ are contained in $P$ and hence $a_1\in\partial\cH_2^+$ and $a_2\in\partial\cH_1^+$.  In particular,
  $a_1\in\cH_2^-$.  Therefore, the closed tropical halfspace $(a_1-a_2)+\cH_2^+$ with apex~$a_1$ which is parallel to
  $\cH_2^+$ is a closed tropical halfspace containing~$P$.  Since $\cH_1^+$ is minimal,
  $\cH_1^+\subseteq(a_1-a_2)+\cH_2^+$ and hence $\cH_1^+\subseteq\cH_2^+$.  Symmetrically, $\cH_2^+\subseteq\cH_1^+$,
  and the claim follows.
\end{proof}

\begin{rem}\label{rem:om}
  The familiarity of the names for the objects defined could inspire the question whether tropical polytopes and, more
  generally, point configurations in the tropical projective space can be studied in the framework of oriented matroids.
  However, as the example in Figure~\ref{fig:full-triangle} shows, the boundaries of tropical halfspaces spanned by a
  given set of points do \emph{not} form a pseudo-hyperplane arrangement, in general.
\end{rem}

\section{Tropical Determinants and Their Signs}

For algorithmic approaches to ordinary polytopes it is crucial that the incidence of a point with an affine hyperplane
can be characterized by the vanishing of a certain determinant expression.  Moreover, by evaluating the sign of that
same determinant, it is possible to distinguish between the two open affine halfspaces which jointly form the complement
of the given affine hyperplane.  This section is about the tropical analog.

Let $M=(m_{ij})\in\RR^{(d+1)\times(d+1)}$ be a matrix.  Then the \emph{tropical determinant} is defined as \[\tdet
M=\bigoplus_{\sigma\in\Sym(d+1)}\bigodot_{i=0}^d m_{i,\sigma(i)} =
\min\SetOf{m_{0,\sigma(0)}+\dots+m_{d,\sigma(d)}}{\sigma\in\Sym(d+1)}.\] Now $M$ is \emph{tropically singular} if the
minimum is attained at least twice, otherwise it is \emph{tropically regular}.  Tropical regularity coincides with the
strong regularity of a matrix studied by Butkovi\v{c}~\cite{MR95a:15025}; see also Burkard and Butkovi\v{c}~\cite{MR2021063}.

The following is proved in Richter-Gebert et al.~\cite[Lemma~5.1]{FirstSteps}.

\begin{prop}\label{prop:singular}
  The matrix $M$ is tropically singular if and only if the $d+1$ points in $\tropPG^d$ corresponding to the rows
  of~$M$ are contained in a tropical hyperplane.
\end{prop}

From the definition of tropical singularity it is immediate that $M$ is tropically regular if and only if its
transpose $M^\tr$ is.  Hence the above proposition also applies to the columns of~$M$.

The \emph{tropical sign} of $\tdet M$, denoted as $\tsgn M$, is either $0$ or $\pm 1$, and it is defined as follows.  If
$M$ is singular, then $\tsgn M=0$.  If $M$ is regular, then there is a unique $\sigma\in\Sym(d+1)$ such that
$m_{0,\sigma(0)}+\dots+m_{d,\sigma(d)}=\tdet M$.  We let the tropical sign of $M$ be the sign of this
permutation~$\sigma$.  See also \cite[\S3.5.1]{MR94b:93001} and the Remark~\ref{rem:bartau} below.

As it turns out the tropical sign shares some key properties with the (sign of the) ordinary determinant.

\begin{prop} Let $M\in\RR^{(d+1)\times(d+1)}$.\label{prop:elementary}
  \begin{enumerate}
  \item If $M$ contains twice the same row (column), then $\tsgn M=0$.
  \item If the matrix $M'$ is obtained from $M$ by exchanging two rows (columns), we have $\tsgn M'=-\tsgn M$.
  \item $\tsgn M^\tr=\tsgn M$.
  \end{enumerate}
\end{prop}

\begin{proof}
  The first property follows from Proposition~\ref{prop:singular}.  The second one is immediate from the definition of
  the tropical sign. And since permuting the rows of a matrix is the same as permuting the columns with the inverse, we
  conclude that $\tsgn M^\tr=\tsgn M$.
\end{proof}

While the behavior of the sign of the ordinary determinant under scaling a row (column) by $\lambda\in\RR$ depends on
the sign of~$\lambda$, the tropical sign is invariant under this operation.  Given $v_0,\dots,v_d\in\RR^{d+1}$ we write
$(v_0,\dots,v_d)$ for the $(d+1)\times(d+1)$-matrix with row vectors $v_0,\dots,v_d$.

\begin{lem}\label{lem:sign-scaling}
  For $v_0,\dots,v_d\in\RR^{d+1}$ and $\lambda_0,\dots,\lambda_d\in\RR$ we have $\tsgn(\lambda_0\odot
  v_0,\dots,\lambda_d\odot v_d)=\tsgn(v_0,\dots,v_d)$.
\end{lem}

\begin{figure}
  \begin{center}
    \subfigure[Non-degenerate case for $\tau=\tau_{(1,0,0),(0,1,0)}$.\label{fig:tau}]{%
      \begin{overpic}[width=.4\textwidth]{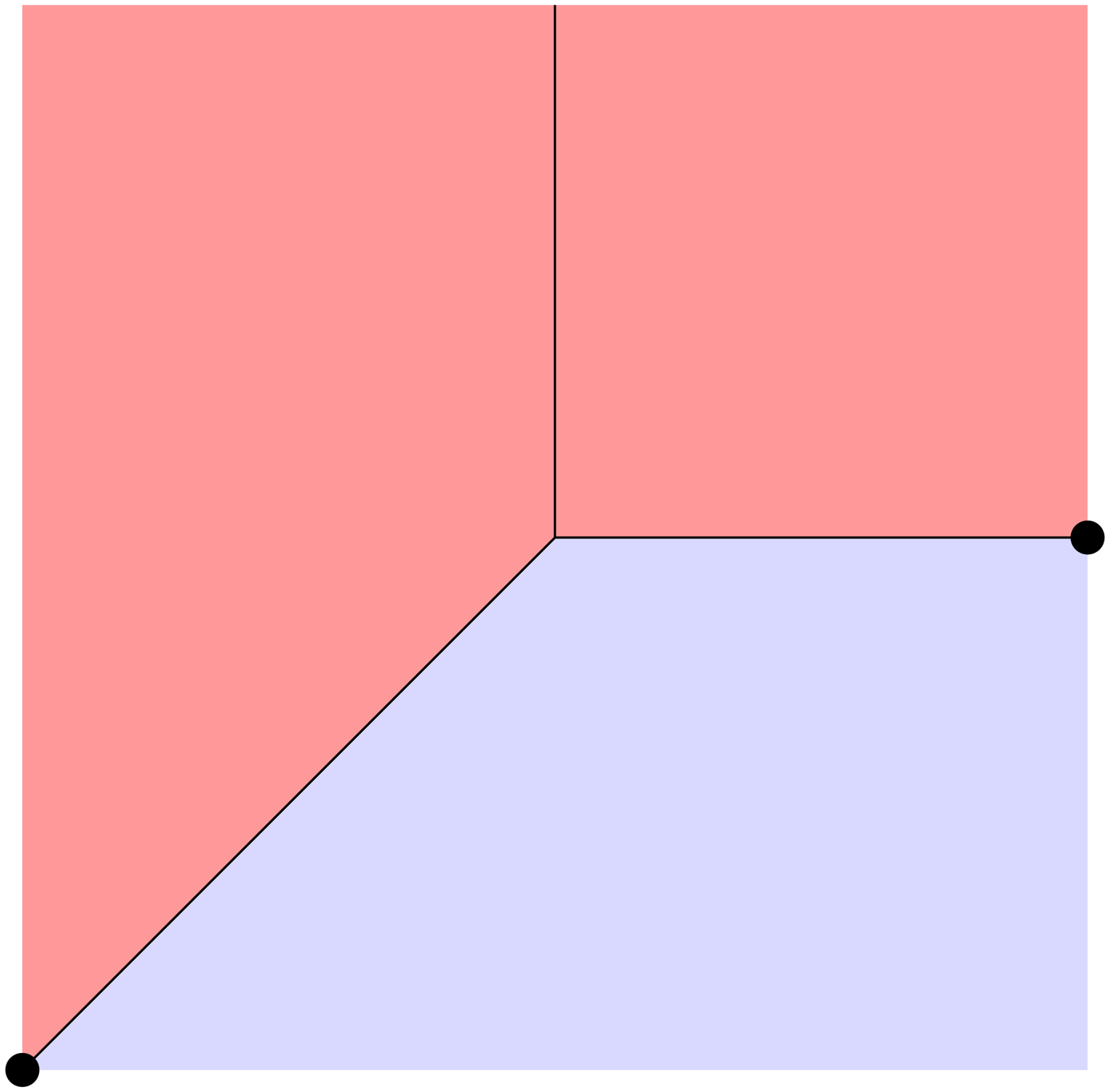}
        \put(-13,15){$(1,0,0)$}
        \put(170,80){$(0,1,0)$}
        \put(20,100){$\tau=1$}
        \put(130,130){$\tau=1$}
        \put(80,30){$\tau=-1$}
      \end{overpic}}
    \quad
     \subfigure[Degenerate case for $\tau'=\tau_{(0,0,0),(0,1,0)}$.\label{fig:tau-dash}]{%
      \begin{overpic}[width=.4\textwidth]{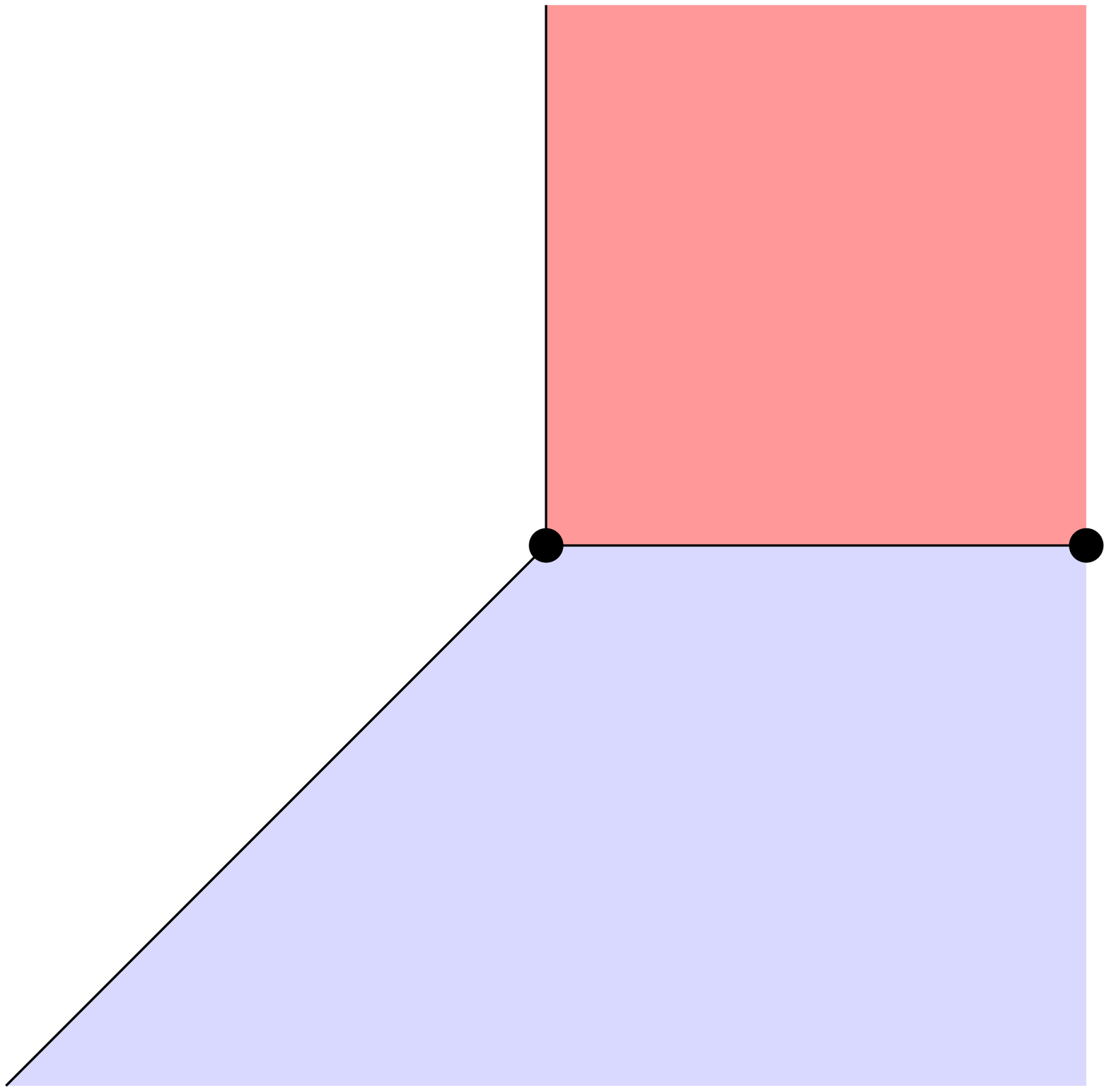}
        \put(75,80){$(0,0,0)$}
        \put(170,80){$(0,1,0)$}
        \put(20,100){$\tau'=0$}
        \put(130,130){$\tau'=1$}
        \put(80,30){$\tau'=-1$}
      \end{overpic}}
    \caption{Values of $\tau_{p,q}$ in $\tropPG^2$ for two different pairs of points.  On the tropical line spanned by
      the two black points the values are zero in both cases.}
  \end{center}
\end{figure}

In fact, $\tsgn$ is a function on $(d+1)$-tuples of points in the tropical projective space~$\tropPG^d$.  For given
$p_1,\dots,p_d$, consider the function \[\tau_{p_1,\dots,p_d}:\tropPG^d\to\{-1,0,1\}:x\mapsto\tsgn(x,p_1,\dots,p_d).\]
Note that we do allow the case where the points $p_1,\dots,p_d$ are not in \emph{general position}, that is, they may be
contained in more than one tropical hyperplane; see the example in Figure~\ref{fig:tau-dash}.

\begin{exmp}\label{exmp:standard-tsgn}
  Consider the real $(d+1)\times(d+1)$-matrix formed of the vertices $-e_0,\dots,-e_d$ of the tropical standard
  simplex~$\Delta^d$.  Then we have \[\tdet(-e_0,\dots,-e_d)=-d,\] and the matrix is regular: The unique minimum is
  attained for the identity permutation, hence \[\tsgn(-e_0,\dots,-e_d)=1,\] or equivalently,
  $\tau_{-e_1,\dots,-e_d}(-e_0)=1$.
\end{exmp}

\begin{prop}
  The function $\tau_{p_1,\dots,p_d}$ is constant on each connected component of the set
  $\tropPG\setminus\tau_{p_1,\dots,p_d}^{-1}(0)$.
\end{prop}

\begin{proof}
  Equip the set $\{-1,0,1\}$ with the discrete topology.  Away from zero the function $\tau_{p_1,\dots,p_d}$ is
  continuous, and the result follows.
\end{proof}

Throughout the following we keep a fixed sequence of points $p_1,\dots,p_d$, and we write $\pi_{ij}$ for the $j$-th
canonical coordinate of~$p_i$.  We frequently abbreviate $\tau=\tau_{p_1,\dots,p_d}$ as well as
$p(\sigma)=\pi_{1,\sigma(1)}+\dots+\pi_{d,\sigma(d)}$ for $\sigma\in\Sym(d+1)$. 

\begin{rem}
  The points $p_1,\dots,p_d$ are in general position if and only if no $d\times d$-minor of the $d\times(d+1)$-matrix
  with entries $\pi_{ij}$ is tropically singular; see \cite[Theorem 5.3]{FirstSteps}.  In the terminology
  of~\cite{TropRank} this is equivalent to saying that the matrix $(\pi_{ij})$ has maximal tropical rank~$d$.
\end{rem}

\begin{thm}\label{thm:halfspace}
  The set $\SetOf{x\in\tropPG^d}{\tau(x)=1}$ is either empty or the union of at most $d$ open sectors of a fixed
  tropical hyperplane.  Conversely, each such union of open sectors arises in this way.
\end{thm}

\begin{proof}
  We can assume that $\tau(x)=1$ for some $x\in\tropPG^d$, since otherwise there is nothing left to prove.  From
  Proposition~\ref{prop:singular} we know that the $d+1$ points $x,p_1,\dots,p_d$ are not contained in a tropical
  hyperplane, and hence they are the vertices of a full tropical $d$-simplex $\Delta=\tconv\{x,p_1,\dots,p_d\}$.
  Consider the facet $F=\tconv\{p_1,\dots,p_d\}$, and let $\cH^+$ be the unique corresponding closed tropical halfspace
  which is minimal with respect to~$\Delta$ and for which we have $\partial\cH^+\cap\Delta=F$.  Let $a$ be the apex
  of~$\cH^+$, and let $a+S_k$ be the open sector containing~$x$.  By construction $a+S_k\subset\cH^+$.
  
  Assume that $\tau(y)\ne\tau(x)$ for some $y\in a+S_k$.  Then there exists a point $z\in[x,y]$ with $\tau(z)=0$.  Since
  sectors are tropically convex, $z\in a+S_k$.  By Proposition~\ref{prop:singular} there exists a tropical
  hyperplane~$\cK$ which contains the points $z,p_1,\dots,p_d$.  Let $\cK^+$ be the minimal closed tropical halfspace of
  the tropical hyperplane~$\cK$ containing $x,p_1,\dots,p_d$.  As $\cH^+$ and $\cK^+$ are both minimal with respect to
  the tropical simplex~$\Delta$, the Proposition~\ref{prop:unique} says that $\cH^+=\cK^+$, and in particular, $a+S_k\ni
  z$ is an open sector of~$\cK$.  The latter contradicts, however, $z\in\cK$.
  
  For the converse, it surely suffices to consider the tropical hyperplane $\cZ$ corresponding to the zero tropical
  linear form, since otherwise we can translate.  We have to prove that for each set $K\subset[d+1]$ with $1\le\#K\le d$
  there is a set of points $u_1,\dots,u_d\in\cZ$ such that
  \[
  \SetOf{x\in\tropPG^d}{\tau_{u_1,\dots,u_d}(x)=1}\ =\ \bigcup\SetOf{S_k}{k\in K}.
  \]
  More specifically, we will even show that, for arbitrary $K$, those $d$ points can be chosen among the
  $\binom{d+1}{2}$ vertices of the second tropical hypersimplex $\Delta_2^d$; see Example~\ref{exmp:Ed}.  Since the
  symmetric group $\Sym(d+1)$ acts on the set of open sectors of~$\cZ$ as well as on the set $\Vert(\Delta_2^d)$, it
  suffices to consider one set of sectors for each possible cardinality $1,\dots,d$.  Let us first consider the case
  where $d$ is odd and $K=\{0,2,4,\dots,d-1\}$ is the set of even indices, which has cardinality $(d+1)/2$.  We set
  \[q_i=\begin{cases} -e_i-e_{i+1} & \text{for $i<d$}\\ -e_0-e_d & \text{for $i=d$}.\end{cases}\]
  If we want to evaluate $\tau_{q_1,\dots,q_d}(x)$ for some point $x\in\tropPG^d\setminus\cZ$ with canonical coordinates
  $(\xi_0,\dots,\xi_d)$, we compute the tropical determinant and the tropical sign of $\tdet(x,q_1,\dots,q_d)$, which in
  canonical row coordinates looks as follows:
  \[
  Q_d=\begin{pmatrix}
    \xi_0  & \xi_1  & \xi_2  & \xi_3  & \xi_4  & \cdots & \xi_d\\
    1      & 0      & 0      & 1      & 1      & \cdots & 1\\
    1      & 1      & 0      & 0      & 1      & \cdots & 1\\
    \vdots & \vdots & \ddots & \ddots & \ddots & \ddots & \vdots\\
    1      & 1      & \cdots & 1      & 0      & 0      & 1\\
    1      & 1      & \cdots & 1      & 1      & 0      & 0\\
    0      & 1      & \cdots & 1      & 1      & 1      & 0
  \end{pmatrix}.
  \]
  Since $x\not\in\cZ$, there is a unique permutation $\sigma_x\in\Sym(d+1)$ such that
  $\tdet(Q)=\xi_{\sigma_x(0)}+q(\sigma_x)$.  We can verify that $\tdet(Q)=0$ in all cases and that
  \[\sigma_x=\begin{cases} (0) & \text{if $x\in S_0$,}\\ (0\ k\ k+1\ \cdots\ d) & \text{if $x\in S_{k}$ for
      $k>0$.}\end{cases}\]  Here we make use of the common cycle notation for permutations, and $(0)$ denotes the
  identity.  For $k>0$ this means that $\sigma_x$ is a $(d+2-k)$-cycle, which is an even permutation if and only if $k$
  is even, since $d$ is odd.  We infer that $\tau_{q_1,\dots,q_d}(x)=1$ if and only if $x\in S_k$ for $k$ even.
  
  We now discuss the case where $\#K\ge(d+1)/2$ and $d$ is arbitrary.  As in the case above, by symmetry, we can assume
  that $K=\{0,2,4,\dots,2(l-1),2l-1,2l,\dots,d\}$ for some $l<\lfloor{d/2}\rfloor$.  We define \[q_i'=-e_0-e_i,\]
  for
  all $i\ge 2l+1$, and we are concerned with the matrix $(x,q_1,\dots,q_l,q_{l+1}',\dots,q_d')$, which, in canonical row
  coordinates, looks like this:
  \[Q_d^l=\begin{pmatrix}
    \xi_0  & \xi_1  & \xi_2  & \xi_3  & \xi_4  & \cdots & \xi_{2l-1} & \xi_{2l} & \cdots & \xi_d\\
    1      & 0      & 0      & 1      & 1      & \cdots & 1          & 1        & \cdots  & 1\\
    1      & 1      & 0      & 0      & 1      & \cdots & 1          & 1        & \cdots  & 1\\
    \vdots & \vdots & \ddots & \ddots & \ddots & \ddots & \vdots     & \vdots   & \ddots & \vdots \\
    1      & 1      & \cdots & 1      & 0      & 0      & 1          & 1        & \cdots & 1\\
    1      & 1      & \cdots & 1      & 1      & 0      & 0          & 1        & \cdots  & 1\\
    0      & 1      & \cdots & 1      & 1      & 1      & 0          & 1        & \cdots & 1\\
    \vdots & \vdots & \ddots & \vdots & \vdots & \vdots & \ddots     & \ddots   & \ddots & \vdots\\
    0      & 1      & \cdots & 1      & 1      & 1      & \cdots     & 1        & 0      & 1\\
    0      & 1      & \cdots & 1      & 1      & 1      & \cdots     & 1        & 1      & 0
  \end{pmatrix}.
  \]
  Note that the upper left $2l\times2l$-submatrix is exactly $Q_l$.  Hence the same reasoning now yields
  \[\sigma_x=\begin{cases} (0) & \text{if $x\in S_0$,}\\ (0\ k\ k+1\ \cdots\ 2l-1) & \text{if $x\in S_{k}$ for
      $k>0$,}\end{cases}\] and $\sigma_x$ is an even permutation if and only if $k$ is even or $k>2l$.
  
  Scrutinizing the matrices $Q_d$ and $Q_d^l$ yields that none of their $d\times(d+1)$-submatrices consisting of all
  rows but the first contains a tropically singular minor.  Equivalently, the points $q_1,\dots,q_d$ as well as the
  points $q_1,\dots,q_l,q_{l+1}',\dots,q_d'$ are in general position.  But then the set
  \[\SetOf{x\in\tropPG^d}{\tau_{q_1,\dots,q_l,q_{l+1}',\dots,q_d'}(x)=-1}\]
  is just the union of the sectors in the
  complement $[d+1]\setminus K$, and since further, according to Proposition~\ref{prop:elementary},
  $\tau_{q_1,\dots,q_d}=-\tau_{q_2,q_1,q_3,\dots,q_l,q_{l+1}',\dots,q_d'}$, the argument given so far already covers the
  remaining case of $\#K<(d+1)/2$.  This completes the proof.
\end{proof}

Now, for the fixed set of points $p_1,\dots,p_d$, we can glue together the connected components of $\tropPG^d\setminus
\tau_{p_1,\dots,p_d}^{-1}(0)$ into two (if $\tau_{p_1,\dots,p_d}\not\equiv0$) large chunks according to their tropical
sign: To this end we define the \emph{closure} of the function $\tau_{p_1,\dots,p_d}$ as follows. Let
$\bar\tau_{p_1,\dots,p_d}(x)=\epsilon$ if there is a neighborhood $U$ of $x$ such that $\tau_{p_1,\dots,p_d}$ restricted
to $U\setminus\tau_{p_1,\dots,p_d}^{-1}(0)$ is identically $\epsilon$; otherwise let $\bar\tau_{p_1,\dots,p_d}(x)=0$.
Clearly, if $\tau_{p_1,\dots,p_d}(x)\ne0$ then $\bar\tau_{p_1,\dots,p_d}(x)=\tau_{p_1,\dots,p_d}(x)$.

Theorem~\ref{thm:halfspace} then implies the following.

\begin{cor}
  The set $\SetOf{x\in\tropPG^d}{\bar\tau(x)=1}$ is empty or a closed tropical halfspace.  Conversely, each closed
  tropical halfspace arises in this way.
\end{cor}

\begin{rem}\label{rem:bartau}
  One can show that $\bar\tau_{p_1,\dots,p_d}(x)=1$ if and only if all optimal permutations, that is, all
  $\sigma\in\Sym(d+1)$ with $\tdet(x,p_1,\dots,p_d)=\xi_{\sigma(0)}+p(\sigma)$ are even.  In this sense our function
  $\bar\tau$ captures the sign of the determinant in the symmetrized min-plus-algebra as defined
  in~\cite[\S3.5.1]{MR94b:93001}.
\end{rem}

\begin{cor}
  For each point $x=(\xi_0,\dots,\xi_d)\in\tropPG^d$ with $\bar\tau_{p_1,\dots,p_d}(x)=0$ there are two permutations
  $\sigma$ and $\sigma'$ of opposite sign such that
  $\tdet(x,p_1,\dots,p_d)=\xi_{\sigma(0)}+p(\sigma)=\xi_{\sigma'(0)}+p(\sigma')$.
\end{cor}

\section{Convex Hull Algorithms in $\tropPG^2$}\label{sec:2d}

For points in the ordinary Euclidean plane the known algorithms can be phrased easily in terms of sign evaluations of
certain determinants.  It turns out that the results of the previous sections can be used to ``tropify'' many ordinary
convex hull algorithms.

In this section we do not use canonical coordinates for points in the tropical projective space, but rather we normalize
by setting the first coordinate to zero.  This way the description of the algorithms can be made in the ordinary affine
geometry language more easily.  In particular, a point in $(0,\xi_1,\xi_2)\in\tropPG^2$ is determined by its
\emph{$x$-coordinate} $\xi_1$ and its \emph{$y$-coordinate} $\xi_2$.  We hope that this helps to see the strong
relationship between the ordinary convex hull problem in $\RR^2$ and the tropical convex hull problem in $\tropPG^2$.
Moreover, this way it may be easier to interpret the illustrations.

Consider a set $S=\{p_1,\dots,p_n\}\subset\tropPG^2$.  Let $\lr(S)$ be the lowest point (least $y$-coordinate) of~$S$
with ties broken by taking the rightmost (highest $x$-coordinate) one.  Similarly, let $\rh(S)$ be the rightmost one
with ties broken by taking the highest, $\hl(S)$ the highest with ties broken by taking the rightmost, and $\lh(S)$ the
leftmost one with ties broken by taking the highest.  Clearly, some of the four points defined may coincide.
If a set of points is in \emph{general position}, that is, for any two points of the input their three rays are pairwise
distinct, then there are unique points with minimum and maximum $x$- and $y$-coordinate respectively.  In this case
there are no ties.

\begin{figure}[htbp]
  \begin{center}
    \begin{overpic}[width=.6\textwidth]{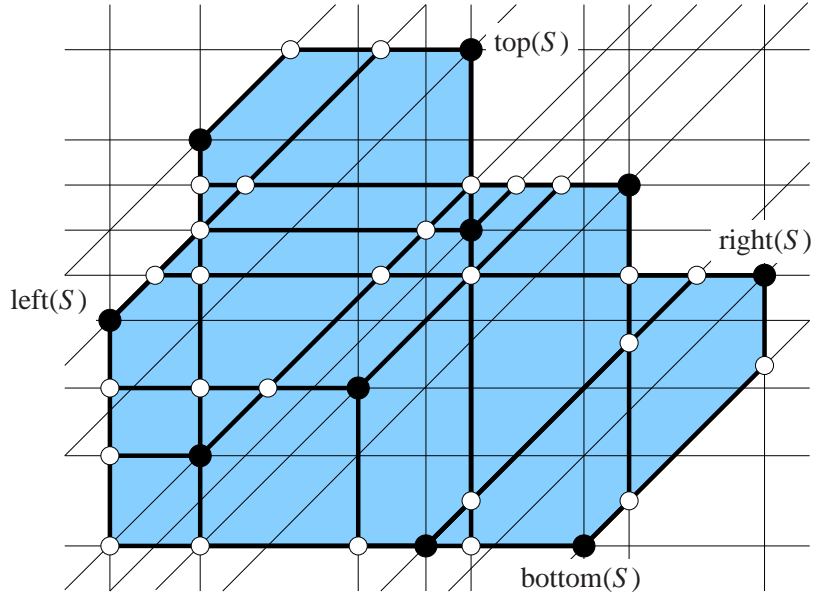}
      \put(170,2){\colorbox{white}{$\lr(S)$}}
      \put(245,130){\colorbox{white}{$\rh(S)$}}
      \put(160,205){\colorbox{white}{$\hl(S)$}}
      \put(-23,107){\colorbox{white}{$\lh(S)$}}
    \end{overpic}
    \caption{Standard affine line arrangement generated by a set of points $S\subset\tropPG^2$, displayed in black.  The
      white points are the pseudovertices on tropical line segments between any two points.  Additionally, the tropical
      convex hull is marked.\label{fig:arrangement}}
  \end{center}
\end{figure}

\begin{figure}[hbtp]
  \begin{center}
    \begin{overpic}[width=.6\textwidth]{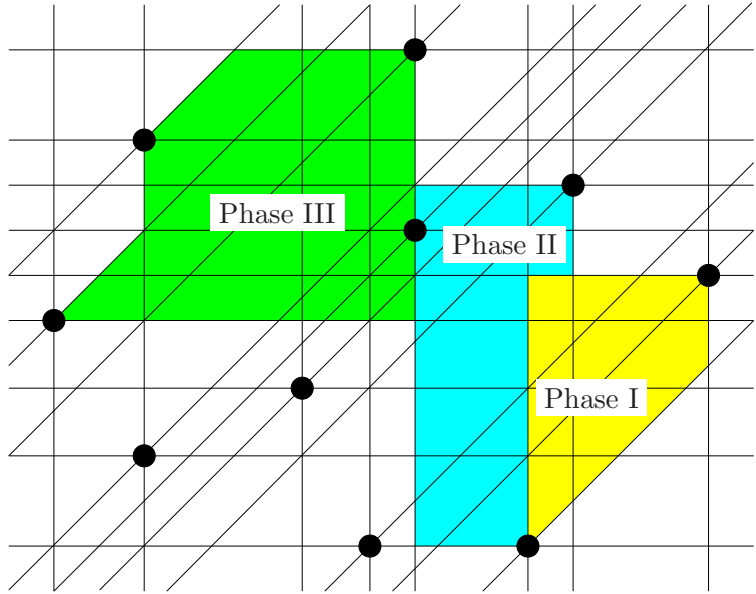}
      \put(200,70){\colorbox{white}{Phase I}}
      \put(165,128){\colorbox{white}{Phase II}}
      \put(77,140){\colorbox{white}{Phase III}}
    \end{overpic}
    \caption{Three phases of the Algorithm~\ref{alg:triple}.\label{fig:phases}}
  \end{center}
\end{figure}

\begin{lem}\label{lem:4vertices}
  The points $\lr(S),\rh(S),\hl(S),\lh(S)$ are vertices of the tropical polygon $\tconv(S)$.  Moreover, $[\lr(S),\lh(S)]$
  is a facet provided that $\lr(S)\ne\lh(S)$.
\end{lem}

\begin{proof}
  By definition, the closed sector $\lr(S)+\bar S_1$ does not contain any point of~$S$ other than $\lr(S)$.  This
  certifies that, indeed, the point $\lr(S)$ is a vertex because of Propositions \ref{prop:vertices}
  and~\ref{prop:zero}.  We omit the proofs of the remaining statements, which are similar.
\end{proof}

Note that due to the special shape of the tropical lines, it is important how to break the ties.  For instance, if two
points have the same lowest $y$-coordinate, then the left one is also on the boundary, but not necessarily a vertex; see
Figure~\ref{fig:arrangement}.

A key difference between tropical versus ordinary polytopes is that the former only have few directions for their
(half-)edges.  This can be exploited to produce convex hull algorithms which do not have a directly corresponding
ordinary version.

Through each point $p=(0,\xi,\eta)\in\tropPG^2$ there is a unique tropical line consisting of three ordinary rays emanating
from~$p$: We respectively call the sets $\SetOf{(0,\xi+\lambda,\eta)}{\lambda\ge0}$, $\SetOf{(0,\xi,\eta+\lambda)}{\lambda\ge0}$,
and $\SetOf{(0,\xi-\lambda,\eta-\lambda)}{\lambda\ge0}$ the \emph{horizontal}, \emph{vertical}, and \emph{skew} ray
through~$x$.  If we have a second point $p'=(0,\xi',\eta')$ then we can compare them according to the relative positions
of their three rays.  This way there is a natural notion of left and right, but there are two notions of above and
below, which we wish to distinguish carefully: $p'$ is \emph{$y$-above} $p$ if $\eta'>\eta$, and it is
\emph{skew-above} if $\eta'-\xi'>\eta-\xi$.

The introduction of the sign of the tropical determinant now clearly allows to take most ordinary $2$-dimensional convex
hull algorithms and produce a ``tropified'' version with little effort.  For instance a suitable tropical version of
Graham's scan provides a worst-case optimal $O(n\log n)$-algorithm.  We omit the details since we describe a different
algorithm with the same complexity.  The commonly expected output of an ordinary convex hull algorithm in two dimensions
is the list of vertices in counter-clockwise order.  As the results in the previous section imply that the combinatorics
of tropical polygons in $\tropPG^2$ does not differ from the ordinary, we adopt this output strategy.

The data structures for the Algorithm~\ref{alg:triple} below is are three doubly-linked lists $L,Y,B$ such that each
input point occurs exactly once in each list.  It is important that all three lists can be accessed at their front and
back with constant cost.

In order to obtain a concise description we assume that the input set $S$ is in general position.  For input not in
general position the notions of left, right, above, and below have to be adapted as above.  The complexity of the
algorithm remains the same.

\begin{algorithm}[htbp]
  \caption{Triple sorting algorithm.\label{alg:triple}}
  \dontprintsemicolon
  \Input{$S\subset\tropPG^2$ finite}
  \Output{list of vertices of $\tconv(S)$ in counter-clockwise order}
  \BlankLine
  sort $S$ from left to right and store the result in list~$L$\;
  sort $S$ from $y$-below to $y$-above and store the result in list~$Y$\;
  sort $S$ from skew-below to skew-above and store the result in list~$B$\;
  \BlankLine
  $H\leftarrow\text{empty list}$; $v\leftarrow\text{front}(Y)$; $w\leftarrow\text{next}(v,Y)$\;
  \While{$w$ $y$-below {\rm back}$(L)$}{
    \If{$w$ skew-below $v$}{
      $v\leftarrow w$; append $v$ to $H$
    }
    $w\leftarrow\text{next}(w,Y)$
  }
  $v\leftarrow\text{back}(L)$; append $v$ to $H$\;
  \BlankLine
  $w\leftarrow\text{previous}(v,L)$\;
  \While{$w$ to the right of {\rm back}$(Y)$}{
    \If{$w$ $y$-above $v$}{
      $v\leftarrow w$;  append $v$ to $H$
    }
    $w\leftarrow\text{previous}(w,L)$
  }
  $v\leftarrow\text{back}(Y)$; append $v$ to $H$;
  \BlankLine
  $w\leftarrow\text{back}(B)$\;
  \While{$w$ skew-above {\rm front}$(L)$}{
    \If{$w$ to the left of $v$}{
      $v\leftarrow w$; append $v$ to $H$
    }
    $w\leftarrow\text{previous}(w,B)$
  }
  $v\leftarrow\text{front}(L)$; append $v$ to $H$\;
  \BlankLine
  \If{$v\ne{\rm front}(Y)$}{
    append front$(Y)$ to $H$
  }
  \Return $H$\;
\end{algorithm}

\begin{prop}
  The Algorithm~\ref{alg:triple} correctly computes the vertices of the tropical convex hull in counter-clockwise
  order.
\end{prop}

\begin{proof}
  The algorithm has an initialization and three phases, where each phase corresponds to one of the three while-loops;
  for an illustration see Figure~\ref{fig:phases}.  In the first phase all the vertices between $\lr(S)$ and
  $\rh(S)$ are enumerated, in the second phase the vertices between $\rh(S)$ and $\hl(S)$, and in the third phase the
  vertices between $\hl(S)$ and $\lh(S)$.
  
  By Lemma~\ref{lem:4vertices} the point~$\text{front}(Y)=\lr(S)$ is a vertex.  Throughout the algorithm the following
  invariant in maintained: $v$ is the last vertex found, and $w$ is an input point not yet processed, which is a
  candidate for the next vertex in counter-clockwise order.  We have a closer look at Phase~I, the remaining being
  similar.  If $w$ is a vertex between $\lr(S)$ and $\rh(S)$ then it will be $y$-above of~$v$, hence we process the
  points according to their order in the sorted list~$Y$.  However, none of those vertices can be $y$-above
  $\text{back}(L)=\rh(S)$, therefore the stop condition.  Under these conditions $w$ is a vertex if and only if $w$ is
  skew-below the tropical line segment $[v,\rh(S)]$.
\end{proof}

The worst-case complexity of the algorithm based on triple sorting is $O(n\log n)$.  If, however, the points are
uniformly distributed at random, say, in the unit square, then by applying Bucket Sort, we can sort the input in an
expected number of $O(n)$ steps; see Cormen et al.~\cite{MR2002e:68001}.

If only few of the input points are actually vertices of the convex hull, then it is easy to beat an $O(n\log n)$
algorithm.  For ordinary planar convex hulls the Jarvis' march algorithm is known as an easy-to-describe method which is
output-sensitive in this sense.  We sketch a ``tropified'' version, we will be instrumental later.  Its complexity is
$O(nh)$, where $h$ is the number of vertices.

\begin{algorithm}[H]
  \caption{Tropical Jarvis' march algorithm.}
  \dontprintsemicolon
  \Input{$S\subset\tropPG^2$ finite}
  \Output{list of vertices of $\tconv(S)$ in counter-clockwise order}
  \BlankLine
  $v_0\leftarrow\lr(S)$; $v\leftarrow v_0$; $H\leftarrow\text{empty list}$\;
  \Repeat{$v=v_0$}{
    $w\leftarrow\text{some point in~$S$}$\;
    \For{$p\in S\setminus\{w\}$}{
      \If{$\bar\tau_{v,w}(p)=-1$ or ($\bar\tau_{v,w}(p)=0$ and $\norm{p-v}>\norm{w-v}$)}{
        $w\leftarrow p$\;
      }
    }
    $v\leftarrow w$; $S\leftarrow S\setminus\{v\}$; append $v$ to $H$\;
  }
  \Return $H$\;
\end{algorithm}

In the ordinary case Chan~\cite{MR97e:68133} gave a worst-case optimal $O(n\log h)$-algorithm, which is based on a
combination of Jarvis' march and a divide-and-conquer approach.  We sketch how the same ideas can be used to obtain an
$O(n\log h)$ convex hull algorithm in~$\tropPG^2$.  If we split the input into $\lceil n/m\rceil$ parts of size at most
$m$, then we can use our $O(n\log n)$ algorithm based on triple sorting to compute the $\lceil n/m\rceil$ hulls in
$O((n/m)(m\log m))=O(n\log m)$ time.  Now we use Jarvis' march to combine the $\lceil n/m\rceil$ tropical convex hulls
into one.  The crucial observation is that each vertex of the big tropical polygon is also a vertex of one of the
$\lceil n/m\rceil$ small tropical polygons.  Therefore, in order to compute the next vertex of the big tropical polygon
in the counter clockwise order, we first compute the tropical tangent through the current vertex to each of the small
tropical polygons.  In each small tropical polygon this can be done by binary searching the vertices in their cyclic
order; this requires $O(\log m)$ steps per small tropical polygon and per vertex of the big tropical polygon.  Summing
up this gives a total of $O(n\log m+h((n/m)\log m))=O(n(1+h/m)\log m)$ operations.  That is to say, if we could know
the number of vertices of the big tropical polygon beforehand, then we could split the input into portions of size at
most~$h$, thus arriving at a complexity bound of $O(n\log h)$.  But this can be achieved by repeated guessing as has
been suggested by Chazelle and Matou\v{s}ek~\cite{MR95k:68229}.

We summarize our findings in the following result, which is identical to the ordinary case.  Note that, as in the
ordinary case, one has an $\Omega(n\log n)$ worst case lower bound for the complexity of the two-dimensional tropical
convex hull problem which comes from sorting.  In this sense our output-sensitive algorithm is optimal.

\begin{thm}
  The complexity of the problem to compute the tropical convex hull of $n$ points in $\tropPG^2$ with $h$ vertices is as
  follows:

  \begin{enumerate}
  \item There is an output-sensitive $O(n\log h)$-algorithm.
  \item There is an algorithm which requires expected linear time for random input.
  \end{enumerate}
\end{thm}

\section{Concluding Remarks}

One of the main messages of this paper is that, with suitably chosen definitions, it is possible to build up a theory of
tropical polytopes quite similar to the one for ordinary convex polytopes.  But, of course, very many items are still
missing.  We list a few open questions, and the reader will easily find more.

\begin{enumerate}
\item How are the face lattices of tropical polytopes related to the face lattices of ordinary convex polytopes?  In
  particular, how do the face lattices of tropical polytopes in $\tropPG^3$ look alike?
\item It is known~\cite[Lemma~22]{TropConvex} that the tropical convex hull of $n$ points in $\tropPG^d$ is the bounded
  subcomplex of some $(n+d)$-dimensional unbounded ordinary convex polyhedron (defined in terms in inequalities).  Hence
  the tropical convex hull problem can be reduced to solving a (dual) ordinary convex hull problem, followed by a search
  of the bounded faces in the face lattice.  The question is: How does an intrinsic tropical convex hull algorithm look
  alike that works in arbitrary dimensions?  Here \emph{intrinsic} means that the algorithm should not take a detour via
  that face lattice computation in the realm of ordinary convex polytopes.  While the complexity status of the ordinary
  convex hull problem is notoriously open (output-sensitive with varying dimension) it is well conceivable that the
  tropical analog is, in fact, easier.  An indication may be the easy to check certificate in
  Proposition~\ref{prop:zero} which leads to a simple and fast algorithm for discarding the non-vertices among the input
  points, a task which is polynomially solvable in the ordinary case, but which requires an LP-type oracle.
\item What is the proper definition of a tropical triangulation?  Such a definition would say that a tropical
  triangulation should be a subdivision into tropical simplices which meet properly.  The precise notion of meeting is
  subtle, however.  While it is obvious that the standard intersection as sets does not do any good, the more refined
  way by extending the $\sqcap$ operation also leads to surprising examples.  A meaningful definition of a tropical
  triangulation should lead to one solution of the previous problem.
\item Can the dimension of an arbitrary tropical polytope, which is not necessarily full, computed in polynomial time?
  Here \emph{dimension} means the same as \emph{tropical rank}.  In fact, this is Question~(1) at the end of the
  paper~\cite{TropRank}.
\item As mentioned in Remark~\ref{rem:om} point configurations in the tropical projective space do not generate an
  oriented matroid in the usual way.  But does there exist a more general notion than an oriented matroid which
  encompasses the tropical case?
\end{enumerate}

\section*{Acknowledgments}

I am grateful to Julian Pfeifle and Francisco Santos for many helpful conversations on the contents of this paper and,
especially, for suggesting better names for some of the notions defined than I could have come up with.  I am indebted
to Bernd Sturmfels for requiring the correction of a few minor errors.  Finally, I would like to thank the anonymous
referee who helped alleviate my ignorance of previously published results on max-plus-algebras.

\bibliographystyle{plain}
\bibliography{main}
\end{document}